\documentclass[12pt]{article}
\usepackage[english]{babel}
\usepackage{amsmath}
\usepackage{amssymb, amsthm, amsfonts}

\newtheorem{theorem}{\hskip 0.5 cm Theorem}

\newtheorem{lemma}{\hskip 0.5 cm Lemma}

\newtheorem{note}{\hskip 0.5 cm Note}

\newtheorem{statement}{\hskip 0.5 cm Proposition}

\newenvironment{pr}{\par\noindent\textbf{Proof.}}{\hfill$\square$}

\begin{document}
\title{Isomonodromic deformation of the linear differential system in the Birkhoff standard form}
\author{Yulia Bibilo\footnote{Department of Theory of Information Transmission and Control, 
Institute for Information Transmission Problems, Russian Academy of Sciences,
Bolshoy Karetny per. 19, Moscow, 127994, Russia. Email: y.bibilo@gmail.com}} 
\date{}
\maketitle
\begin{abstract}
We consider a linear meromorphic system in the Birkhoff standard form. The construction of the isomonodromic deformation of it proposed by Bolibruch is discussed. This construction has some special characteristics because of resonant irregular singularity at the infinity. 
\end{abstract}

The theory of isomonodromic deformations for Fuchsian systems was developed by L.~Schle\-sin\-ger, L.~Fuchs, R.~Garnier and later by B.~Malgrange, A.~Bolibruch and others. M.~Jimbo, T.~Miwa, K.~Ueno constructed isomonodromic deformations of meromorphic system with Fuchsian and irregular non-resonant singularities $\cite{Jimbo1}$. Malgrange used theory of holomorphic vector bundels with meromorphic connections and improved this result $\cite{Ma1}$ (He constructed isomonodromic deformation with irregular resonant singularities such that each corresponding leading coefficient matrix has no Jordan blocks with equal eigenvalues). Later V.~Heu used Malgrange's approach and constructed isomonodromic deformation for any $(2\times2)$-meromorphic system \cite{Heu}. 
M.~Bertola, M.~Y.~Mo considered case when leading coefficient matrix corresponding to each irregular singularity has no equal Jordan blocks $\cite{Bert_Mo}$.
A.~Bolibruch gave a brief construction  of the isomonodromic deformation of the system with one Fuchsian singularity and one "non-ramified" irregular resonant singularity in $\cite{Bol7}$. The goal of this work is to make some necessary additions and proofs for his construction.

\section{System in the Birkhoff standard form}

Consider a linear differential meromorphic system
\begin{equation}\label{syst0}
\frac{dy}{dz}=A(z)y,\quad A(z)=A_{r-1}z^{r-1}+\ldots+A_1+\frac{A_0}{z},\quad y \in \mathbb{C}^d, \quad z\in \bar{\mathbb{C}},
\end{equation}
here $z=0$ is a Fuchsian singular point, $z=\infty$ is an irregular singular point.

Two systems $$\frac{dy}{dz}=A(z)y,\qquad \frac{d\tilde{y}}{dz}=\tilde{A}(z)\tilde{y}$$ are called locally meromorphically equivalent at $z=\infty$, if
there is meromorphic matrix $\Gamma(z)$, $\det \Gamma(z) \not \equiv 0$ such, that coefficient matrices $A(z)$ and $\tilde{A}(z)$ satisfy equality $$\tilde{A}(z)=\Gamma^{-1}(z)A(z)\Gamma(z)-\Gamma^{-1}(z)\frac{d\Gamma(z)}{dz}$$ for the gauge transformation $$\tilde{y}=\Gamma(z)y.$$

Analogously, two systems  $$\frac{dy}{dz}=A(z)y,\qquad \frac{d\tilde{y}}{dz}=\tilde{A}(z)\tilde{y}$$  are formally equivalent at $z=\infty$, if there is formal Laurent series $\hat{\Gamma}(z)=\sum_{j=-k}^{\infty}\Gamma_j z^{-j}$ such, that $$\tilde{y}=\hat{\Gamma}(z)y$$ transforms $A(z)$ to $\tilde{A}(z)$. 

Poinkare rank $r$ of the system at the infinity is called minimal Poinkare rank, if Poinkare ranks of all locally meromorphically equivalent systems are equal or large than $r$.

The system $(\ref{syst0})$ is said to be in the Birkhoff standard form, if it's Poinkare rank $r$ at the infinity is minimal.

\section{Formal and analytic solutions}

The system $(\ref{syst0})$ has the following formal fundamental matrix at the infinity $\cite{Wasow,Balser1}$:
\begin{equation}\label{ff_matr0}
\hat{Y}_{\infty}(z)=\hat{F}(z) z^L e^{Q^0(z)},
\end{equation}
where $\hat{F}(z)$ is a formal Laurent series in $z^{-1/p}$ with finite principal part, $p\in \mathbb{N}$, $Q^0(z)$ is a diagonal matrix which non-zero elements are polynomials in  $z^{1/p}$ with zero free terms. Moreover, analytic continuation $\tilde{Q}^0(z)$ of the matrix $Q^0(z)$ along loop encircling $\infty$ satisfies the relation 
\begin{equation}
\tilde{Q}^0(z)=R^{-1}Q^0(z)R,
\end{equation}
 here $R$ is a constant permutation matrix.
A formal monodromy matrix of the system $(\ref{syst0})$ at the infinity is equal to $\hat{G}_{\infty}=e^{2\pi \mathbf{i} L}=DR$, where $D$ is a constant diagonal matrix.

There is another factorization of the formal fundamental matrix $\cite{Balser1}$:
\begin{equation*}
\hat{Y}_{\infty}(z)=\hat{U}(z)z^J U e^{Q^0(z)},
\end{equation*} 
where 
$\hat{U}(z)$ and $(\hat{U}(z))^{-1}$ are formal Laurent series in $1/z$ with finite principle part,
block-diagonal matrices $Q^0(z,b)$, $J$, $U$ have agreed block structure, every block $Q_j(z)$ of $Q^0(z)$ has the form
\begin{equation*}
Q_j(z)={\rm diag}(q_j(t)I_{s_j},q_j(t\varsigma_j)I_{s_j},\ldots,q_j(t\varsigma_j^{p_j-1})I_{s_j}),
\end{equation*}
where $q_j$ is a polynomial in $t=z^{1/p_j}$, $\varsigma_j=e^{2\pi \mathbf{i}/p_j}$ for some integer $p_j$ not greater then common multiple of $2,3,\ldots,p$, also, polynomial $q_j$ has no constant term . Every block of $J$ has the form
\begin{equation*}
J_j={\rm diag}(J_{s_j},J_{s_j}+1/p_j I_{s_j},\ldots,J_{s_j}+(p_j-1)/p_j I_{s_j}),
\end{equation*}
and every block $U_j$ of $U$ decomposes into blocks of the form
\begin{equation*}
U_{lk}=[\varsigma_j^{(l-1)(k-1)}I_{s_j}],\quad 1\leq l,k \leq p_j.
\end{equation*}
In this way formal monodromy matrix has form $\hat{G}_{\infty}=U^{-1}e^{2\pi \mathbf{i} J}U$. 

If $p=1$ in $Q^0(z)$, irregular singular point of the system $(\ref{syst0})$ is called non-ramified. 

Formal solutions converge to analytic ones in the neighborhood of the Fuchsian singular point $\cite{Wasow}$. It is useful to have following description of the analytic solutions at Fuchsian points in terms of Levelt basis $\cite{Le}$.
Let $E$ denote a matrix which is defined as a logarithm of the monodromy matrix $G_0$:
\begin{equation}
E=\frac{1}{2\pi \mathbf{i}}\ln G_0, \quad 0 \leq Re \rho_j <1, \quad j=1,\ldots,p,
\end{equation}
where $\rho_1,\ldots, \rho_p$ are eigenvalues of the matrix $E$.
Matrix $S$ exists such, that fundamental matrix $Y_0(z)$ in the neighborhood of the Fuchsian singularity $z=0$ has the form
\begin{equation}\label{levelt_y}
Y_0(z)S=U(z)z^{\Lambda}z^{E'},
\end{equation}
matrix $U(z)$ is holomorphic and holomorphically invertible at zero, $\Lambda$ is a integer-valued diagonal matrix, $E'=S^{-1}ES$ is an upper-triangular matrix.

\section{Monodromy data}

Here we will define monodromy data of the system $(\ref{syst0})$ $\cite{Bol3,Bol7,Bol_Malek_Mit}$. Consider arbitrary point $z_0 \in \mathbb{C}\backslash \{0\}$
and analytic continuation $\tilde{Y}(z)$ of the fundamental matrix $Y(z)$ along the simple loop $\gamma$ encircling one singular point and with the beginning at $z_0$. $Y(z)$ and $\tilde{Y}(z)$ are fundamental matrices of the same linear differential system, then they are connected by the relation $\tilde{Y}(z)G=Y(z)$, where $G$ is a non-degenerate constant matrix. Thus, we define the monodromy representation 
\begin{equation}
\chi:\pi_1(\mathbb{C}\backslash\{0\})\rightarrow GL(p,\mathbb{C}),\quad  \chi:[\gamma]\rightarrow G.
\end{equation}
Representation $\chi$ depends only on homotopic class $[\gamma]$ of the loop $\gamma$ \cite{Bol3}.

Monodromy data of the system $(\ref{syst0})$ are completely defined by the set of Stokes matrices $\{S_1,\ldots,S_{N}\}$ at the infinity and by the monodromy matrix $G_{\infty}$ at the infinity. (i.e. $G_{\infty}$ is a matrix corresponding to analytic continuation of the fundamental matrix around $z=\infty$. Monodromy matrix $G_0$ at zero is equal to inverse matrix to $G_{\infty}$, $G_0=G_{\infty}^{-1}$.)

The Stokes matrices can be defined in following way $\cite{Bol_Malek_Mit}$\footnote{see also $\cite{Balser2}$}. Let $l_1 \prec \ldots \prec l_N$ be the singular rays, it means they are the rays from $\infty$ and are labeled in ascending order with respect to the positive orientation of a circle centered at $\infty$, on which some $e^{q_i-q_j}$ has maximal decay. Formal fundamental solution $\hat{Y}_{\infty}$ is multi-summable along every non-singular ray $l$.
Let $l_i^{+},l_i^{-}$ be two rays such that $l_i^{-} \prec l_i \prec l_i^+$ and $l_i$ is only singular ray of $Q^0(z)$ containing in the oriented sector $[l_i^-,l_i^+]$. Let $Y_i^-$ and $Y_i^+$ denote the sums of $\hat{Y}_{\infty}$ along $l_i^-$ and $l_i^+$ respectively. These solutions are connected by $Y_i^+=Y_i^-S_i$ in the neighborhood of $l_i$, where constant matrix $S_i$ is called Stokes matrix. Stokes matrices satisfy relations:
\begin{equation*}
e^{Q^0(z)}S_{k}e^{-Q^0(z)}\sim I, 
\end{equation*}
\begin{equation*}
S_1 \cdot \ldots \cdot S_N \cdot \hat{G}_{\infty}=G_{\infty}.
\end{equation*}

\section{Parametric Sibuya theorem}
We call monodromy data, consisting of matrices $J$, $Q(z,b)$ ($b\in \mathbb{C}^t$), $G_{\infty}$ and Stokes matrices $\{S_1,\ldots,S_N\}$ corresponding to the respective singular directions $l_1 \prec \ldots \prec l_N$ of $Q(z,b)$, is admissible, if\\
(i) there is a natural number $p$;\\
(ii) $Q(z,b)$ and $J$ have agreed block structure, every block of $Q(z,b)$ has the form
\begin{equation*}
Q_j(z,b)={\rm diag}(q_j(t,b)I_{s_j},q_j(t\varsigma,b)I_{s_j},\ldots,q_j(t\varsigma^{p_j-1},b)I_{s_j}),
\end{equation*}
where $q_j$ is a polynomial in $t=z^{1/p_j}$, $\varsigma_j=e^{2\pi \mathbf{i}/p_j}$ for some integer $p_j$ not greater then common multiple of $2,3,\ldots,p$, also, polynomial $q_j$ has no constant term . Every block of $J$ has the form
\begin{equation*}
J_j={\rm diag}(J_{s_j},J_{s_j}+1/p_j I_{s_j},\ldots,J_{s_j}+(p_j-1)/p_j I_{s_j}).
\end{equation*}
(iii) Singular directions $l_1 \prec \ldots \prec l_N$ of $Q(z,b)$ are constant. Stokes matrices for every $b$ satisfy conditions
\begin{equation*}
e^{Q(z,b)}S_{k}e^{-Q(z,b)}\sim I, \quad z\rightarrow \infty, \quad z \in [l_i-\epsilon,l_i+\epsilon],
\end{equation*}
for small positive $\epsilon$,
\begin{equation*}
S_1 \cdot \ldots \cdot S_N \cdot \hat{G}_{\infty}=G_{\infty}.
\end{equation*}

\begin{theorem}\label{th_sib_param}$\cite{Sibuya,Bol_Malek_Mit}$
Let the monodromy data $J$, $Q(z,b)$, $G_{\infty}$, $\{S_1,\ldots,S_N\}$ is admissible,
diagonal elements of the matrix $Q(z,b)$ are polynomials in $z$ with degrees to be equal or less then $r$, and which coefficients are holomorphic in $b$ in $\mathbb{C}^m$. Then there is a neighborhood $D(b^k)\subset \mathbb{C}^m$ and a system $dy=\omega_k y$ for each point $b^k \in \mathbb{C}^m$ such, that:\\
1. coefficient matrix of the system $dy=\omega_k y$ is holomorphic in $b \in D(b^k)$;\\
2. $z=\infty$ is an irregular singular point of the system $dy=\omega_k y$; \\
3. monodromy data of the system $dy=\omega_k y$ at the infinity is equal to the given set $\{S_1,\ldots,S_N\}$, $J$, $Q(z,b)$.
\end{theorem}
Theorem $\ref{th_sib_param}$ in non-parametric case was formulated in $\cite{Sibuya}$ but Stockes matrices were defined in different way. It is also easily follows from the sufficient conditions for generalized Riemann--Hilbert problem positive solution.

\begin{theorem}\label{equivalent_syst}\label{th_Sib2}$\cite{Sibuya}$
Let two meromorphic equations $$\frac{dy}{dz}=A(z)y,\qquad \frac{d\tilde{y}}{dz}=\tilde{A}(z)\tilde{y}$$ are formally equivalent at the infinity, then they are meromorphically equivalent if and only if their Stockes structures at the infinity are equal.
\end{theorem}

\section{ Parametric Savage lemma}

Let $H(K)$ be the space of holomorphic inside $K$ and continuous on $K$ matrix-valued functions, and $H^0(K)$ be the subspace of $H(K)$, which contains only holomorphically invertible inside $K$ matrix-valued functions. Denote also $O_0=\{z\in \mathbb{C}:|z|<R\}$, $O_{\infty}=\{z\in \bar{\mathbb{C}}:|z|>r\}$, $R>r>0$, $K=O_0 \cap O_{\infty}$.

\begin{lemma}\label{trivial_bund}$\cite{Bol3}$
Let $F(z,b)$ belongs to $H^0(K \times D(b^0))$ and the inequality  $\left\| F-I  \right\| < \epsilon$ is true for quite small $\epsilon >0$, then matrix-valued functions $W(z,b)\in H^0(O_0\times D(b^0))$, $U(z,b)\in H^0(O_{\infty} \times D(b^0))$ exist such, that $F(z,b)=W^{-1}(z,b)U(z,b)$.
\end{lemma}

\begin{note}\label{note_cont}$\cite{Bol3}$
Let $F(z,b) \in H^0(K \times D_{\delta_0}(b^0))$.
For any $\epsilon>0$ there is $\delta >0$ such, that $||F^{-1}(z,b^0)F(z,b)-I||<\epsilon$ in $K \times D_{\delta}(b^0)$.
\end{note}


\begin{lemma}\label{Sovage}$\cite{Bol3}$
Let $\hat{F}(z)$ is a formal Laurent series at the infinity with finite principle part. Then there are a formal Laurent series $\hat{U}(z)$, $\det \hat{U}(\infty) \neq 0$, a rational matrix-valued function $\Gamma(z)\in H^0(\mathbb{C})$ and an integer matrix $M$ such, that $\Gamma(z) \hat{F}(z)=z^M \hat{U}(z)$. 
\end{lemma}

\begin{lemma}
Let conditions of the theorem $\ref{th_sib_param}$ are fulfilled and exponential part $Q(z,b)$ of the formal fundamental solution is polynomial in $z$ ($p=1$). Then the formal fundamental solution $\hat{Y}_{\infty}$ at the infinity of the system $dy=\omega_k y$, defined by the theorem $\ref{th_sib_param}$, has the form $\hat{X}^k(z,b)=\hat{U}^k(z,b)z^L e^{Q(z,b)}$, where $\hat{U}^k(z,b)$ is a formally invertible Taylor series in $1/z$ and analytic in $b$.
\end{lemma}\label{lem_sib_ad}
\begin{pr}
There is system $dy=\vartheta_k y$ defined by the theorem $\ref{th_sib_param}$. It has fundamental solution $Y_{\infty}(z,b)=F(z,b)z^Le^{Q(z,b)}$ at the infinity, where $F(z,b)$ is a meromorphic function in $z^{-1}$ and analytic function in $b$. According to the note $\ref{note_cont}$ $||F^{-1}(z,b^0)F(z,b)-I||$ is small. Then $$F(z,b)=F(z,b^0)W^{-1}(z,b)U(z,b),$$ where $U(z,a)\in H^0(O_{\infty} \times D(b^0))$ and $\Gamma(z,b)=F(z,b^0)W^{-1}(z,b)$ is meromorphic in $z$ and analytic in $b$. Let apply gauge transformation $\tilde{y}=\Gamma(z,b)$ to the system $dy=\vartheta_k y$, we will get new system $dy=\omega_k y$ with the fundamental matrix $\tilde{Y}_{\infty}(z,b)=U(z,b)z^Le^{Q(z,b)}$ at the infinity. Then formal fundamental matrix has the form $\hat{X}^k(z,b)=\hat{U}^k(z,b)z^L e^{Q(z,b)}$, where $\hat{U}^k(z,b)$ is a formally invertible Taylor series in $1/z$ and analytic in $b$, because of $U(z,b)\sim \hat{U}^k(z,b)$, $z\rightarrow \infty$ in some sector.
\end{pr}

\section{Malgrange theorem}

\begin{theorem}\label{Malgrange_th}$\cite{Ma1,Bol3,Bol7}$
Let $M(z,b)$ be a holomorphically invertible matrix function in $K \times T$, where $T$ is a connected analytic manifold and $K=O_0 \cap O_{\infty}$. Suppose that for some $t^0 \in T$ one has $M(z,t^0)=(W_0(z))^{-1}U_0(z)$ with $W_0(z)\in H^0(O_{\infty})$, $W_0(\infty)=I$, $U_0(z) \in H^0(O_0)$. Then there exist an analytic subset $\Theta \subset T$ of codimension one and unique holomorphic mappings
\begin{equation*}
W:O_{\infty} \times (T\backslash \Theta) \rightarrow GL(p,\mathbb{C}),
U:O_{0} \times (T\backslash \Theta) \rightarrow GL(p,\mathbb{C}),
\end{equation*}
such that the following conditions are fulfilled\\
(i) $W,W^{-1}$ are meromorphic along $O_{\infty} \times \Theta$ and $W(\infty,t)\equiv I$;\\
(ii) $U,U^{-1}$ are meromorphic along $O_0 \times \Theta$;\\
(iii) $M(z,t)=(U(z,t))^{-1} W(z,t)$.
\end{theorem}

\section{Isomonodromic family constructing}
A.~Bolibruch proposed construction of isomonodromic deformation of $(\ref{syst0})$ in $\cite{Bol7}$, one can find it in the proof of the theorem $\ref{isom_deform_constr}$ below. The goal of this work is to make some necessary additions.

Let matrix-valued function $Q(z,b)$ is given such that it satisfies the conditions:\\
(*) $Q(z,b)$ is a diagonal matrix and all diagonal elements of $Q(z,b)$ are polynomials in $z^{\frac{1}{p}}$ with coefficients are analytic in $b\in \mathbb{C}^m$ and free terms are zeroes;\\
(**) singular directions of $Q(z,b)$ are constant. 

The problem is to construct isomonodromic family in the Birkhoff standard form such that conditions are hold:\\
1. for any $b$ the formal fundamental matrix has exponential part given by $Q(z,b)$;\\
2. monodromy data, i.e. monodromy matrix $G_{\infty}$ and Stokes matrices are constant;\\
3. given system $(\ref{syst0})$ is included in the isomonodromic family and $Q(z,0)=Q^0(z)$.

This problem has solution in the case of non-ramified irregular singular point.

\begin{statement}\label{stat_1}
Let system $(\ref{syst0})$ has non-ramified irregular singular point, $Q(z,b)$ satisfies $(*),(**)$ and $Q(z,b)$ is polynomial in $z$ ($p=1$). Then
\begin{enumerate}
\item There is a neighborhood $O_{\infty}$ of the infinity such that following statements are hold.
\item For any $b^k\in \mathbb{C}^m$ there exist a small disk $D(b^k)$ and isomonodromic family $dy=\omega^k y$ over $O_{\infty}\times D(b^k)$ such that it has an irregular rank $r$ singularity at the infinity, it's exponential part of the formal solution is equal to $Q(z,b)$, it's monodromy data coincide with monodromy data of $(\ref{syst0})$.
\item If $D(b^k) \cap D(b^m) \neq  \varnothing$, then formal fundamental matrices $\hat{Y}^k(z,b)$, $\hat{Y}^m(z,b)$ of $dy=\omega^k y$, $dy=\omega^m y$ satisfy condition
\begin{equation}\label{F_cond}
\hat{Y}^k(z,b)(\hat{Y}^k(z,b))^{-1} \sim C, \quad z \rightarrow \infty,
\end{equation}
where $C$ is a constant in $z$ and non-degenerate matrix.
\item When $b=0$ given system $(\ref{syst0})$ is equal to the system $dy=\omega^0 y$ over $O_{\infty}\times D(0)$ (index $0$ corresponds to disk $D(0)$).
\end{enumerate}
\end{statement}

\begin{pr} 
Let $l_1 \prec \ldots \prec l_N$ be singular directions of $Q(z,b)$ and $\{S_1,\ldots,S_N\}$ be corresponding Stokes matrices of $(\ref{syst0})$, $J,G_{\infty}$ be monodromy data of the system $(\ref{syst0})$.   

According to the Sibuya parametric theorem (theorem $\ref{th_sib_param}$) there is a neighborhood $O_{\infty}$ of the infinity such that for any $b^k\in \mathbb{C}^m$ there is a small disk $D(b^k)$ and an isomonodromic family
\begin{equation}\label{local_isom}
dy=\vartheta^k y,\quad \vartheta^k=A(z,b)dz=\left(A_r(b)z^{r-1}+\ldots+\frac{A_0(b)}{z}+\frac{A_{-1}(b)}{z^2}+\ldots\right)dz,
\end{equation}
such that the exponential part of the formal fundamental solution matrix is equal to $Q(z,b)$ and $(\ref{local_isom})$ has monodromy data $\{S_1,\ldots,S_N\}$, $J,G_{\infty}$. Formal fundamental matrix $\hat{X}^k(z,b)$ of $(\ref{local_isom})$ has form $\hat{X}^k(z,b)=\hat{U}^k(z,b)z^L e^{Q(z,b)}$, where $\hat{U}^k(z,b)$ is a formally invertible Taylor series in $z$ and analytic in $b$ (by lemma $\ref{lem_sib_ad}$).

We modify the system $dy=\vartheta^0 y$ over $O_{\infty} \times D(0)$ by a gauge transformation so that it includes given system $(\ref{syst0})$ when $b=0$. According to the Savage lemma $\ref{Sovage}$ formal series $\hat{F}(z)$ of $(\ref{ff_matr0})$ can be written in the form $\hat{F}(z)=\hat{U}(z) \Gamma^{-1}(z) z^M$, $\hat{U}(z)$ is a formal Laurent series, $\det \hat{U}(\infty) \neq 0$,  $\Gamma(z)\in H^0(\mathbb{C})$ is a rational matrix-valued function, and $M$ is an integer matrix. Then for any $D(b^k)$ consider family with formal fundamental matrix 
$$\hat{Y}^k(z,b)=\hat{F}^k(z,b) z^L e^{Q(z,b)}, \quad \hat{F}^k(z,b)=\hat{U}^{k}(z,b) \Gamma^{-1}(z) z^M$$ 
and let denote it by
\begin{equation}\label{x_isom}
dy=\omega^k y.
\end{equation}
The system $(\ref{x_isom})$ is analytic in $b \in D(b^k)$ and it has the same monodromy data. Indeed, by the Savage lemma $\ref{Sovage}$, for each fixed $b\in D(b^k)$ there is a meromorphic over $O_{\infty}$ matrix $\Gamma_1(z)$ such that $\hat{U}^{k}(z,b) \Gamma^{-1}(z) z^M = \Gamma_1(z)z^{M_1} \hat{U}^k_1(z,b)$, where $M_1$ is a constant diagonal integer matrix, $\hat{U}^k_1(z,b)$ is a formal and formally invertible Taylor series.
Thus systems $dy=\omega^k y$ and $dy=\vartheta^k y$ are meromorphically equivalent and they have the same Stokes matrices by the Sibuya theorem $\ref{th_Sib2}$. So the system $dy=\omega^0 y$ over $O_{\infty} \times D(0)$ is equal to the given system $(\ref{syst0})$ when $b=0$.

Consider two disks $D(b^k)$, $D(b^m)$ such, that $D(b^k) \cap D(b^m) \neq  \varnothing$, and corresponding isomonodromic families $(\ref{x_isom})$. Formal fundamental matrices $\hat{Y}^k(z,b)$, $\hat{Y}^m(z,b)$ of these families are satisfy the condition:
\begin{equation}\label{F_cond}
\hat{Y}^k(z,b)(\hat{Y}^k(z,b))^{-1}=\hat{F}^k(z,b)(\hat{F}^m(z,b))^{-1} \sim C
\end{equation}
for $z \rightarrow \infty$, where $C$ is a constant non-degenerated matrix.
\end{pr}

\begin{theorem}\label{isom_deform_constr}
Let matrix $Q(z,b)$ satisfies conditions (*),(**) and $Q(z,b)$ is polynomial in $z$ ($p=1$).
Then there is an isomonodromic family in the Birkhoff standard form satisfying conditions 1,2,3.
\end{theorem}

\begin{pr}
1. Let $O_{\infty}$, $\{D(b^k)\}$, $\{dy=\omega^k y\}$ are defined by the proposition $\ref{stat_1}$.

Consider covering $\{D(b^k)\}$ of $\mathbb{C}^m$. Suppose $D(b^k) \cap D(b^j) \neq \varnothing$.
Consider sector $S$ being in $O_{\infty}$ with vertex at $z=\infty$ such that it's solution less then $\pi/N$. We can find analytic fundamental matrices $Y^k(z,b)$, $Y^j(z,b)$ of $dy=\omega^k y$, $dy=\omega^j y$ accordingly such, that $Y^k(z,b) \sim \hat{Y}^k(z,b)$ in $S \times D(b^k)$, $Y^j(z,b) \sim \hat{Y}^j(z,b)$ in $S \times D(b^j)$.

Next we consider function $g_{kj}(z,b)=Y^k(z,b) (Y^j(z,b))^{-1} \in H(O_{\infty} \times D(b^k) \cap D(b^j))$, it is holomorphically invertible in $O_{\infty} \times D(b^k) \cap D(b^j)$. Indeed, $g_{kj}(z,b) \in H^0(O_{\infty}\backslash \{\infty \} \times D(b^k) \cap D(b^j))$ as a product of non-degenerate matrices, and $g_{kj}(z,b)$ is holomorphically invertible in $\{\infty \} \times D(b^k) \cap D(b^j)$ because of relation $(\ref{F_cond})$.

In that way, covering $\{O_{\infty}\times D(b^k)\}$ of manifold $O_{\infty}\times \mathbb{C}^m$ and cocycle  $\{g_{kj}\}$ define bundle $P$. The bundle $P$ is holomorphically trivial as a bundle over contractible Stein manifold. Thus, there is bundle trivialization $\{X^k(z,b)\}$ such that $X^0(z,0)=I$ (index $0$ corresponds to disk $D(0)$).

Consider function 
\begin{equation} 
Y(z,b)=(X^k(z,b))^{-1} Y^k(z,b)
\end{equation}
in $O_{\infty} \times \mathbb{C}^m$, it is analytic outside $\{\infty\} \times \mathbb{C}^m$. Moreover, it is isomonodromic fundamental matrix of the family $dy=\omega_{\infty} y$, $\omega_{\infty}=dY(z,b)Y^{-1}(z,b)$ in $O_{\infty}\times \mathbb{C}^m$. Note, initial system is included to the constructed family when $b=0$.

2. Write fundamental matrix $Y(z)$ of the initial system $(\ref{syst0})$ in the form $(\ref{levelt_y})$, i.e.
\begin{equation}
Y(z)=U_0(z)z^{\Lambda} z^{E'}S^{-1},
\end{equation}
where $S$ ia a constant matrix, $E'=S^{-1} E S$ is an upper-triangular matrix, $A$ is a diagonal matrix, with integer diagonal elements creating non-increasing sequence, and $U_0(z)$ is a holomorphically invertible in $\mathbb{C}$ matrix.

Let write matrix $Y(z,b)S$ in similar form in $O_{\infty} \times \mathbb{C}^m$
\begin{equation}
Y(z,b)S=T(z,b)z^Az^{E'}.
\end{equation}
And then we consider holomorphic vector bundle $F$ over $\bar{\mathbb{C}}\times \mathbb{C}^m$ which is defined by covering $\{O_{\infty} \times \mathbb{C}^m,\mathbb{C} \times \mathbb{C}^m\}$ 
and cocicle $g_{\infty0}=T(z,b)$. Differential forms $\omega_{\infty}$, $\omega_0=(A+z^AE'z^{-A})z^{-1}dz$ define connection $\nabla$ on the bundle. $\omega_{\infty}$ can be calculated by $\omega_{\infty}=dYY^{-1}=dg_{\infty 0}(g_{\infty 0})^{-1}+g_{\infty 0}\omega_{0} (g_{\infty 0})^{-1}$.

3. Vector bundle $F$ is holomorphically trivial on $\bar{\mathbb{C}}\times \{0\}$, then for matrix-valued function $T(z,b)$ theorem $\ref{Malgrange_th}$ are fulfilled. So, there is analytic subset $\Theta \subset \mathbb{C}^m$ and matrix-valued functions $\Gamma(z,b)$, $U(z,b)$ such, that

(i) $\Gamma(z,b)$ is holomorphically invertible in $O_{\infty} \times (\mathbb{C}^m \backslash \Theta)$ and meromorphic in $O_{\infty} \times \Theta$ and $\Gamma(z,0)=I$.

(ii) $U(z,b)$ is holomorphically invertible in $\mathbb{C} \times (\mathbb{C}^m \backslash \Theta)$ and meromorphic in $\mathbb{C} \times \Theta$.

(iii) $\Gamma(z,b)T(z,b)=U(z,b)$.

Consider matrix-valued function $X(z,b)=\Gamma(z,b)Y(z,b)$. Differential form $$\omega(z,b)=dX(z,b)X^{-1}(z,b)$$ is holomorphic outside zero and infinity. For every fixed $b \in \mathbb{C}^m \backslash \Theta$ differential form is equal to $\omega=A(z,b)dz$, $A(z,b)$ is a polynomial in $z$ with degree $r-1$. Polynomial coefficients $A_{i}(b)$ are holomorphic in $b$ in $\mathbb{C}^m \backslash \Theta$ and meromorphic in $\Theta$. So, given system $(\ref{syst0})$ is included to the isomonodromic family.

\end{pr}

\section*{Acknowledgments}
Y.~Bibilo acknowledges the support of the Russian Foundation for Basic Research (grant no. RFBR 14-01-00346 А and  grant no. RFBR 14-01-31145 mol\underline{ }a).


\begin{thebibliography}{99}
		
		\bibitem{Bol3} 		 
		A.~A. Bolibruch, {\it Inverse Monodromy Problems in the Analytic Theory of Differential Equations}, MCCME, Moscow, 2009 (in Russian).
		
		\bibitem{Bol7}
    A.~Bolibruch, {\it Inverse problems for linear differential equations with meromorphic coefficients},
		CRM Proceeding and Lecture Notes, {\bfseries 31} (2002), pp. 3--25.
	
		\bibitem{Bol_Malek_Mit}		 
		A.~A.~Bolibruch, S.~Malek, C.~Mitschi,		{\it On the generalized Riemann--Hilbert problem with irregular singularities},		
		Expo. Math., {\bfseries 24}:3 (2006), pp. 235--272.	
			
		\bibitem{Balser1}
    W.~Balser, W.~B.~Jurkat, D.~A.~Lutz, {\it A general theory of invariants for meromorphic differential equations. I. Folmal invariants},
		Funkcialaj Ekvacioj, {\bfseries 22}:2 (1979), pp. 197--221.
		
		\bibitem{Balser2}     
		W.~Balser, W.~B.~Jurkat, D.~A.~Lutz,     {\it A general theory of invariants for meromorphic differential equations. II. Proper invariants}, Funkcialaj Ekvacioj, {\bfseries 22}:2 (1979), pp. 257--283.
		
		\bibitem{Bert_Mo}
    M.~Bertola, M.~Y.~Mo, {\it Isomonodromic deformation of resonant rational connections},
		Int. Math. Res. Papers, {\bfseries 2005}:11 (2005), pp. 565--635.
				
		\bibitem{Jimbo1}    
		M.~Jimbo, T.~Miwa, K.~Ueno, {\it Monodromy preserving deformation of linear ordinary differential equations with rational coefficientes}, Physica D, {\bfseries 2} (1981), pp. 306--352.
		
		\bibitem{Heu}
		V.~Heu, {\it Universal isomonodromic deformations of meromorphic rank 2 connections on curves},
		Ann. Inst. Fourier, Grenoble, {\bfseries 60}:2 (2010), pp. 515--549. 
		
		\bibitem{Le}
		A.\,Levelt, {\it Hypergeometric functions. II}, 
		Proc. Konikl. Nederl. Acad. Wetensch. Ser.~A {\bfseries 64} (1961), pp. 373--385.
	  		
		
		\bibitem{Ma1}
		B.\,Malgrange, {\it Sur les d\'eformations isomonodromiques.	I. Singularit\'es r\'eguli\`eres}, 
		Progr. Math. {\bfseries 37} (1983), pp. 401--426.

			
		
		\bibitem{Wasow}
		W.~Wasow, {\it Asymptotic expansions for ordinary differential equations}, Wiley, New~York–London–Sydney, 1965.

		\bibitem{Sibuya}    
		Y.~Sibuya,		{\it Stokes phenomena},     
		Bulletin of the American Mathematical Society, {\bfseries 83}:5 (1977), pp. 1075--1077. 
		
\end{thebibliography}
\end{document}